\documentclass[11pt]{amsart}

\usepackage[margin=1in]{geometry}
\usepackage{amsmath,amssymb,mathtools}
\usepackage{booktabs}
\usepackage{microtype}
\usepackage{titlesec}
\usepackage{tikz}
\usetikzlibrary{arrows.meta,calc,decorations.pathreplacing,positioning}
\usepackage[colorlinks=true,linkcolor=blue,citecolor=blue,urlcolor=blue,
  pdftitle={Block and Residue Exclusion in Kernel-Root Expansions for Two-Barrier Walks},
  pdfauthor={Cetin Hakimoglu-Brown},
  pdfkeywords={bounded lattice paths, absorbing barriers, kernel method, Schur functions, root-of-unity filter}]{hyperref}
\usepackage[nameinlink,noabbrev]{cleveref}

\titleformat{\section}
  {\normalfont\large\bfseries\raggedright}
  {\thesection.}{0.65em}{}
\titlespacing*{\section}{0pt}{2.2ex plus 0.6ex minus 0.2ex}{1.0ex}
\titleformat{\subsection}
  {\normalfont\normalsize\bfseries\raggedright}
  {\thesubsection.}{0.65em}{}
\titlespacing*{\subsection}{0pt}{1.8ex plus 0.5ex minus 0.2ex}{0.8ex}

\makeatletter
\def\@settitle{\begin{center}%
  \baselineskip18pt\relax
  \Large\bfseries
  \@title
  \end{center}%
}
\makeatother

\newtheorem{theorem}{Theorem}[section]
\newtheorem{proposition}[theorem]{Proposition}
\newtheorem{lemma}[theorem]{Lemma}
\newtheorem{corollary}[theorem]{Corollary}
\theoremstyle{remark}
\newtheorem{remark}[theorem]{Remark}

\newcommand{\Q}{\mathbb{Q}}
\newcommand{\cS}{\mathcal{S}}
\newcommand{\cL}{\mathcal{L}}
\newcommand{\ord}{\operatorname{ord}}
\newcommand{\Z}{\mathbb{Z}}

\title[Block and residue exclusion]{Block and Residue Exclusion in Kernel-Root Expansions for Two-Barrier Walks}
\author{Cetin Hakimoglu-Brown}
\address{Berkeley, California, USA}
\email{mathemails@proton.me}
\date{}
\subjclass[2020]{Primary 05A15; Secondary 05E05, 60G50}
\keywords{directed lattice paths, absorbing barriers, bounded walks, kernel method, Schur functions, Puiseux series, root-of-unity filter}

\begin{document}

\begin{abstract}
We consider the two-barrier walk on the integers with steps $+y$ and $-b$, stopped on entering either absorbing set. Its finite-state generating functions are rational and admit the standard Schur-function representation in the roots of the normalized kernel $P_z(X)=X^{b+y}-2X^b+z$.

Our contribution is a two-stage pruning method for constructing these root expansions without computing terms that cannot affect the final polynomial. First, the Weyl expansion is grouped into complete small-root/large-root blocks, and a $z$-adic valuation bound removes every block beginning above the denominator degree. Second, surviving blocks are split into residue subblocks and projected onto complete cyclic root-label orbits; the resulting congruences eliminate nontrivial characters without any termwise rationality test. In the $(+3,-2)$, $m=12$ example only one mixed-block degree survives, while for $(+5,-2)$, $m=22$, twenty of twenty-one Weyl terms are excluded before expansion.
\end{abstract}

\maketitle
\markboth{Cetin Hakimoglu-Brown}{Block and Residue Exclusion}

\section{Introduction}

Directed paths confined to a finite interval, or equivalently to a slit between two parallel lines, have rational generating functions.  The transfer-matrix description is immediate, while the kernel method expresses the same functions in terms of the roots of a characteristic polynomial.  Bousquet-M\'elou developed rectangular Schur denominators for bounded excursions~\cite{BousquetMelou}; Banderier and Nicod\`eme treated walks below one wall and between two walls~\cite{BanderierNicodeme}; and Khalid and Prellberg obtained a general skew Schur representation for directed paths with arbitrary starting and ending heights in a slit~\cite{KhalidPrellberg}.  The relation with banded Toeplitz determinants is discussed by Alexandersson~\cite{Alexandersson}.

The small/large-root form of the kernel method is classical; see Banderier and Flajolet~\cite{BanderierFlajolet}.  The general finite-slit generating function is therefore known.  The issue studied here is different: \emph{how much of a kernel-root expansion must actually be computed?}  A root formula may be a sum of many infinite algebraic series even though the final numerator and denominator are finite polynomials.  Direct expansion consequently produces terms that are later removed by cancellation.

For a two-step model, the kernel roots separate near the origin into small roots, which have positive Puiseux valuation, and large roots, which tend to nonzero constants after a monomial normalization.  This suggests two complementary filters.  First, Weyl terms can be grouped into complete small-root/large-root blocks and compared with the finite degree window of the determinant.  Second, the Lagrange expansions inside a surviving block can be separated by residue classes and summed over complete cyclic root-label orbits.

An earlier version of this manuscript used these ideas in a less invariant form and, in particular, imposed a termwise rationality condition on individual trigonometric phases.  The present version replaces that condition by complete character-orbit sums and makes the finite two-barrier model, the role of the parameter $m$, and the connection with the transfer matrix explicit from the outset.

The combined construction is
\[
\boxed{
\text{block valuation}
\;\longrightarrow\;
\text{residue congruences}
\;\longrightarrow\;
\text{complete orbit sum}
\;\longrightarrow\;
\text{finite truncation}.}
\]
The two examples at the end show the two possible outcomes.  In the $(+3,-2)$ case with $m=12$, the mixed block meets the denominator at exactly one degree.  In the $(+5,-2)$ case with $m=22$, no mixed block reaches the polynomial range, and only the single all-large Weyl term is needed.

\section{The walk, the barriers, and the parameter \texorpdfstring{$m$}{m}}\label{sec:model}

Fix coprime integers
\[
    1\le b<y,
\]
and an integer $m>b$.  Let $(S_k)_{k\ge0}$ be the discrete-time walk
\begin{equation}\label{eq:walk}
    S_{k+1}=\begin{cases}
        S_k+y,&\text{with probability }1/2,\\
        S_k-b,&\text{with probability }1/2.
    \end{cases}
\end{equation}
The finite two-barrier process is obtained by declaring the sets
\begin{equation}\label{eq:absorbing-sets}
    \mathcal L=\{0,1,\ldots,b-1\},
    \qquad
    \mathcal U_m=\{m,m+1,\ldots,m+y-1\}
\end{equation}
absorbing.  The transient state set is therefore
\begin{equation}\label{eq:transient-set}
    \mathcal T_m=\{b,b+1,\ldots,m-1\}.
\end{equation}
Thus $m$ is \emph{not} the number of transient states: it is the first state in the upper absorbing set.  The number of transient states is
\begin{equation}\label{eq:transient-count}
    |\mathcal T_m|=m-b.
\end{equation}
The extra $b-1$ states on the left and $y-1$ states on the right record the possible overshoots produced by steps of size $b$ and $y$.

For a boundary payoff $g$ on $\mathcal L\cup\mathcal U_m$, let
\[
    \tau=\inf\{k\ge0:S_k\in\mathcal L\cup\mathcal U_m\}
\]
and define the ordinary hitting-time generating function
\begin{equation}\label{eq:physical-gf}
    F_s(t)=\mathbb E_s\!\left[t^\tau g(S_\tau)\right],
    \qquad s\in\mathcal T_m.
\end{equation}
For example, to count absorption at the inner left state $b-1$, take
\[
    g(b-1)=1,
    \qquad
    g(j)=0\quad(j\ne b-1).
\]
The examples in \cref{sec:example12,sec:example22} use this payoff and start from
\begin{equation}\label{eq:startstate}
    s=m-1,
\end{equation}
the uppermost transient state.

Let $Q_m$ be the $(m-b)\times(m-b)$ transition matrix restricted to $\mathcal T_m$, and let $r_g$ be the vector of one-step transition probabilities from transient states to absorbing states, weighted by $g$.  Then
\begin{equation}\label{eq:resolvent}
    \mathbf F(t)=t(I-tQ_m)^{-1}r_g.
\end{equation}
In particular, every $F_s(t)$ is rational and has denominator dividing $\det(I-tQ_m)$.

The characteristic equation for the interior recurrence
\[
    F_s(t)=\frac t2\bigl(F_{s+y}(t)+F_{s-b}(t)\bigr)
\]
is
\begin{equation}\label{eq:physical-kernel}
    K_t(\lambda)=t\lambda^{b+y}-2\lambda^b+t=0.
\end{equation}
For the root analysis it is convenient to make the monomial change
\begin{equation}\label{eq:normalization}
    X=t^{1/y}\lambda,
    \qquad
    z=t^{(b+y)/y}.
\end{equation}
Then
\begin{equation}\label{eq:kernel-normalization}
    K_t(t^{-1/y}X)
    =t^{-b/y}\bigl(X^{b+y}-2X^b+z\bigr).
\end{equation}
We shall therefore work with the normalized kernel
\[
    P_z(X)=X^{b+y}-2X^b+z.
\]
Although $z$ is a convenient Puiseux normalization rather than the ordinary time variable, the denominator below contains only powers of $z^y$.  Hence the substitution
\begin{equation}\label{eq:zytn}
    z^y=t^{b+y}
\end{equation}
returns an ordinary polynomial in $t$.

\begin{figure}[htbp]
\centering
\vspace{0.6em}
\begin{tikzpicture}[font=\small,>=Latex]
  \draw[fill=gray!15,rounded corners=1pt] (0,0) rectangle (2.6,1.15);
  \draw[rounded corners=1pt] (2.6,0) rectangle (9.0,1.15);
  \draw[fill=gray!15,rounded corners=1pt] (9.0,0) rectangle (12.3,1.15);

  \node at (1.30,0.68) {$0,\ldots,b-1$};
  \node at (5.80,0.68) {$b,\ldots,m-1$};
  \node at (10.65,0.68) {$m,\ldots,m+y-1$};

  \node[below=6pt] at (1.30,0) {lower absorbing set};
  \node[below=6pt] at (5.80,0) {transient states};
  \node[below=6pt] at (10.65,0) {upper absorbing set};

  \draw[decorate,decoration={brace,amplitude=5pt}]
      (2.68,1.43)--(8.92,1.43)
      node[midway,above=8pt] {$m-b$ transient states};

  \node[below=28pt] at (2.60,0) {first transient state $b$};
  \node[below=28pt] at (9.00,0) {first upper absorbing state $m$};
\end{tikzpicture}
\vspace{0.9em}
\caption{State-space convention. The parameter $m$ is the first state in the upper absorbing set, not the number of transient states. A transient position $j$ moves to $j-b$ or $j+y$, and the walk stops as soon as either shaded region is reached.}
\label{fig:state-space}
\vspace{0.8em}
\end{figure}

\begin{figure}[htbp]
\centering
\vspace{1.4em}
\begin{tikzpicture}[font=\small,>=Latex,x=1.18cm,y=0.38cm]
  \fill[gray!15] (-0.20,-0.35) rectangle (5.75,1.45);
  \fill[gray!15] (-0.20,11.65) rectangle (5.75,14.65);

  \draw[->] (0,0)--(6.05,0) node[right] {time $k$};
  \draw[->] (0,0)--(0,15.15) node[above] {position $S_k$};

  \node at (2.25,0.62) {lower absorbing set $\{0,1\}$};
  \node at (3.55,13.55) {upper absorbing set $\{12,13,14\}$};
  \draw[dashed] (0,12)--(5.72,12);
  \node[anchor=north east] at (5.65,11.88) {barrier begins at $m=12$};

  \draw[very thick,-{Latex[length=2.4mm]}]
    plot coordinates {(0,11)(1,9)(2,7)(3,5)(4,3)(5,1)};
  \foreach \x/\yy in {0/11,1/9,2/7,3/5,4/3,5/1}
    {\fill (\x,\yy) circle (2.25pt);}

  \node[anchor=east] at (-0.15,11.0) {start $s=m-1=11$};
  \node[anchor=south east] at (5.02,1.35) {$S_5=1$};
  \node[below=10pt] at (2.55,0) {five consecutive steps of size $-2$};

  \draw[dashed,thick,-{Latex[length=2.4mm]}] (0.10,11.10)--(1.10,14.10);
  \fill (1.10,14.10) circle (2.25pt);
  \node[anchor=west] at (1.38,14.32) {$+3$ from $11$ gives upper absorption};
\end{tikzpicture}
\vspace{1.3em}
\caption{A sample path for the $(+3,-2)$ walk with $m=12$. Starting from $s=11$, five consecutive $-2$ steps reach the lower absorbing state $1$. The dashed alternative shows the one-step jump from $11$ into the upper absorbing set.}
\label{fig:model}
\vspace{1.2em}
\end{figure}
\section{The two-barrier denominator}

Put
\[
    n=b+y.
\]
With the model and normalization of \cref{sec:model}, we use the normalized kernel
\begin{equation}\label{eq:kernel}
    P_z(X)=X^n-2X^b+z.
\end{equation}
This is precisely the normalized kernel obtained from the physical recurrence in \eqref{eq:kernel-normalization}.  We keep the variable $z$ throughout the root analysis and return to the ordinary time variable $t$ by \eqref{eq:zytn}.

Let $x_1,\ldots,x_n$ be the roots of \eqref{eq:kernel}.  Define the alternant
\begin{equation}\label{eq:alternant}
 A_m(z)=
 \det\!\begin{pmatrix}
 1&1&\cdots&1\\
 x_1&x_2&\cdots&x_n\\
 \vdots&\vdots&&\vdots\\
 x_1^{b-1}&x_2^{b-1}&\cdots&x_n^{b-1}\\
 x_1^m&x_2^m&\cdots&x_n^m\\
 \vdots&\vdots&&\vdots\\
 x_1^{m+y-1}&x_2^{m+y-1}&\cdots&x_n^{m+y-1}
 \end{pmatrix}
\end{equation}
and put
\begin{equation}\label{eq:Ddef}
    D_m(z)=\frac{A_m(z)}{\Delta(x_1,\ldots,x_n)},
    \qquad
    \Delta(x_1,\ldots,x_n)=\prod_{i<j}(x_j-x_i).
\end{equation}
Up to a nonzero scalar normalization, $D_m$ is the common denominator appearing in the two-barrier absorption generating functions.

\begin{proposition}[Schur form]\label{prop:schur}
One has
\begin{equation}\label{eq:schur}
    D_m(z)=s_{(m-b)^y}(x_1,\ldots,x_n).
\end{equation}
\end{proposition}

\begin{proof}
Reverse the rows in the numerator and denominator of \eqref{eq:Ddef}.  The descending exponent sequence in the numerator is
\[
    m+y-1,m+y-2,\ldots,m,b-1,b-2,\ldots,0.
\]
Subtracting the staircase $n-1,n-2,\ldots,0$ gives the partition $((m-b)^y,0^b)$.  The bialternant formula gives \eqref{eq:schur}; see, for example, \cite[Ch.~I]{Macdonald}.
\end{proof}

The parameter $m$ in \eqref{eq:alternant} is exactly the barrier parameter from \cref{sec:model}: the rectangular partition has width $m-b$, the number of transient states.  The normalized determinant is the Schur form of the finite transfer determinant.

Let $e_j=e_j(x_1,\ldots,x_n)$.  Comparing \eqref{eq:kernel} with the monic polynomial having roots $x_i$ gives
\begin{equation}\label{eq:especial}
 e_j=0\quad(1\leq j\leq n-1,\ j\neq y),
 \qquad
 e_y=(-1)^{y+1}2,
 \qquad
 e_n=(-1)^n z.
\end{equation}
Of course $e_0=1$.  The sparsity of this specialization gives both a support restriction and a degree bound.

\begin{proposition}[Polynomial window]\label{prop:window}
For $\gcd(b,y)=1$,
\begin{equation}\label{eq:window}
    D_m(z)\in\Q[z^y],
    \qquad
    \deg_zD_m\leq
    N_m:=y\left\lfloor\frac{m-b}{b+y}\right\rfloor.
\end{equation}
\end{proposition}

\begin{proof}
The Schur polynomial $s_{(m-b)^y}$ has weighted degree $y(m-b)$ when $e_j$ is assigned weight $j$.  Under \eqref{eq:especial}, a nonzero monomial has the form $e_y^a e_n^k$ with integers $a,k\geq0$, and satisfies
\[
    ya+nk=y(m-b).
\]
Here coprimality is used through
\[
    \gcd(n,y)=\gcd(b+y,y)=\gcd(b,y)=1,
\]
so $y\mid k$.  Write $k=yq$.  Since $a\geq0$,
\[
    nyq\leq y(m-b),
\]
so $q\leq\lfloor(m-b)/n\rfloor$.  As $\deg_z$ of such a monomial is $k=yq$, this gives \eqref{eq:window}.
\end{proof}

\begin{remark}[Sharpness of the window]\label{rem:windowsharp}
The upper bound in \eqref{eq:window} is exact in both worked examples and in further computations for several coprime pairs $(y,b)$.  This suggests the general identity
\[
    \deg_zD_m=y\left\lfloor\frac{m-b}{b+y}\right\rfloor,
\]
but exactness is not needed for any exclusion result below: a certified upper bound is sufficient.  We leave the general nonvanishing of the top admissible coefficient as an open question.
\end{remark}

\begin{proposition}[Transfer-matrix normalization]\label{prop:transferbridge}
Let $L=m-b$.  Under the substitution $z^y=t^{b+y}$,
\begin{equation}\label{eq:transferbridge}
    D_m(z)
    =(-1)^{(y+1)L}2^L\det(I-tQ_m).
\end{equation}
\end{proposition}

\begin{proof}
The matrix $2(I-tQ_m)$ is the $L\times L$ banded Toeplitz matrix with diagonal entries $2$, a $-t$ diagonal at displacement $+y$, and a $-t$ diagonal at displacement $-b$, truncated at the two boundaries.  The standard banded Toeplitz--Schur determinant identity identifies its determinant with the rectangular Schur specialization in \eqref{eq:schur}; see, for example, \cite{Alexandersson}.  The monomial change \eqref{eq:normalization} gives $z^y=t^{b+y}$, and \cref{prop:window} shows that the substitution is single-valued on $D_m$.  Finally, at $t=0$ the transfer determinant is $1$, whereas
\[
    D_m(0)=\bigl((-1)^{y+1}2\bigr)^L,
\]
which fixes the scalar and sign in \eqref{eq:transferbridge}.
\end{proof}

\section{Complete small-root and large-root blocks}

At $z=0$,
\[
    P_0(X)=X^b(X^y-2).
\]
Hensel factorization separates the roots into
\[
    \cS(z)=\{r_1(z),\ldots,r_b(z)\},
    \qquad
    \cL(z)=\{R_1(z),\ldots,R_y(z)\},
\]
where
\begin{equation}\label{eq:rootorders}
    \ord_z r_i=\frac1b,
    \qquad
    \ord_z R_j=0.
\end{equation}
The leading coefficients of the small roots are distinct $b$th roots of $1/2$, so
\begin{equation}\label{eq:smalldifference}
    \ord_z(r_i-r_j)=\frac1b
    \qquad(i\neq j).
\end{equation}
Every small--large and large--large difference is a unit at $z=0$.

The rectangular Weyl formula gives
\begin{equation}\label{eq:weyl}
 D_m(z)=
 \sum_{\substack{I\subseteq\{1,\ldots,n\}\\ |I|=y}}
 \frac{\left(\prod_{i\in I}x_i\right)^m}
 {\displaystyle\prod_{\substack{i\in I\\j\notin I}}(x_i-x_j)}.
\end{equation}
For $0\leq c\leq b$, let $B_{m,c}(z)$ be the sum of the terms in \eqref{eq:weyl} for which $I$ contains exactly $c$ small roots.  Then
\begin{equation}\label{eq:blocksum}
    D_m(z)=\sum_{c=0}^{b}B_{m,c}(z),
\end{equation}
and the $c$th block contains
\begin{equation}\label{eq:blocksize}
    \binom bc\binom yc
\end{equation}
Weyl terms.

Although an individual Weyl term is generally a Puiseux series, a complete block is regular and has cyclotomic support.

\begin{lemma}[Block regularity and support]\label{lem:blockregularity}
For every $0\leq c\leq b$,
\begin{equation}\label{eq:blocksupport}
    B_{m,c}(z)\in\Q[[z^y]].
\end{equation}
\end{lemma}

\begin{proof}
Let $A_{m,c}$ be the corresponding part of the Laplace expansion of $A_m$.  It is alternating separately in the small roots and in the large roots, and hence is divisible by $\Delta(\cS)\Delta(\cL)$.  Since
\[
 \Delta(\cS\cup\cL)
 =\pm\Delta(\cS)\Delta(\cL)
   \prod_{\substack{r\in\cS\\R\in\cL}}(R-r),
\]
and the cross product is a unit, $B_{m,c}$ is a separately symmetric expression in the roots of the two Hensel factors.  Therefore $B_{m,c}\in\Q[[z]]$.

Let $\omega^y=1$.  The identity
\[
    P_{\omega^b z}(\omega X)=\omega^bP_z(X)
\]
shows that multiplication of every root by $\omega$ carries the two root sets for $z$ to the corresponding root sets for $\omega^b z$.  Each Weyl term is homogeneous of total root degree $y(m-b)$, so
\[
    B_{m,c}(\omega^b z)=B_{m,c}(z).
\]
Because $b$ is invertible modulo $y$, $\omega^b$ runs through all $y$th roots of unity.  Hence only powers of $z^y$ occur.
\end{proof}

\begin{theorem}[Block valuation]\label{thm:blockvaluation}
For $0\leq c\leq b$,
\begin{equation}\label{eq:rawvaluation}
    \ord_z B_{m,c}
    \geq\frac{c(m-b+c)}{b}.
\end{equation}
Together with \eqref{eq:blocksupport}, this gives
\begin{equation}\label{eq:roundedvaluation}
    B_{m,c}(z)=O\!\left(z^{L_{m,c}}\right),
    \qquad
    L_{m,c}=y\left\lceil
       \frac{c(m-b+c)}{by}
      \right\rceil.
\end{equation}
\end{theorem}

\begin{proof}
In a Weyl term from the $c$th block, the numerator contains $c$ small roots and has order $cm/b$.  The denominator contains exactly $c(b-c)$ small--small cross differences, each of order $1/b$ by \eqref{eq:smalldifference}; all other cross differences are units.  Thus every term has order
\[
    \frac{cm-c(b-c)}b
    =\frac{c(m-b+c)}b.
\]
Summation cannot lower the valuation.  Rounding upward to the next multiple of $y$ gives \eqref{eq:roundedvaluation} by \cref{lem:blockregularity}.
\end{proof}

\begin{corollary}[Whole-block exclusion]\label{cor:blockexclusion}
If
\begin{equation}\label{eq:blockcriterion}
    L_{m,c}>N_m,
\end{equation}
then $B_{m,c}$ cannot affect the denominator polynomial.  More precisely,
\begin{equation}\label{eq:blocktruncation}
    D_m(z)=
    [z^{\leq N_m}]
    \sum_{\substack{0\leq j\leq b\\j\neq c}}B_{m,j}(z),
\end{equation}
where $[z^{\leq N}]$ denotes truncation after degree $N$.
\end{corollary}

\begin{proof}
The block has no coefficient in the polynomial window $0\leq k\leq N_m$.
\end{proof}

\begin{remark}
An excluded block is generally not identically zero.  Its algebraic tail cancels against the tails of other blocks above the degree of $D_m$.  The point is that this cancellation need not be calculated.
\end{remark}

The terminal block is always excluded.

\begin{corollary}[All-small block]\label{cor:allsmall}
For every $m>b$, the block $B_{m,b}$ lies above the degree window.
\end{corollary}

\begin{proof}
From \eqref{eq:roundedvaluation},
\[
    L_{m,b}=y\left\lceil\frac{m}{y}\right\rceil\geq m,
\]
whereas
\[
    N_m\leq\frac{y(m-b)}{b+y}<m.
\]
\end{proof}

\section{Residue subblocks and complete character orbits}\label{sec:residues}

The valuation test removes entire blocks.  We now refine each surviving block by the residue classes occurring in the Lagrange expansions of the roots.  The first step is to derive explicitly the ``strings'' that occur in this refinement.

\subsection{From a Weyl term to a string}\label{subsec:strings}

Fix a Weyl index set $I$ in \eqref{eq:weyl}, and put
\[
    J=I^{\mathrm c}=\{u_1,\ldots,u_b\}.
\]
Thus the string variables come from the \emph{complement} of the $y$-element Weyl set.  This point fixes the index convention used below.  If $I$ contains $c$ small roots, then $J$ contains
\begin{equation}\label{eq:jgfromc}
    j=b-c\quad\text{small roots},
    \qquad
    g=c\quad\text{large roots},
\end{equation}
and therefore $j+g=b$.

Let
\[
    T_I(z)=
    \frac{\left(\prod_{i\in I}x_i\right)^m}
    {\displaystyle\prod_{i\in I,\,j\in J}(x_i-x_j)}
\]
be the corresponding Weyl term.  Since $P_z'(u)=\prod_{x\neq u}(u-x)$, multiplication over $u\in J$ gives
\begin{equation}\label{eq:crossderivative}
 \prod_{u\in J}P_z'(u)
 =(-1)^{by+\binom b2}
   \Delta(u_1,\ldots,u_b)^2
   \prod_{i\in I,\,j\in J}(x_i-x_j).
\end{equation}
Also $\prod_{i=1}^{n}x_i=(-1)^nz$.  Hence
\begin{equation}\label{eq:weylcomplement}
 T_I(z)
 =\sigma_{b,y,m}\,z^m
  \Delta(u_1,\ldots,u_b)^2
  \prod_{\ell=1}^{b}\frac{u_\ell^{-m}}{P_z'(u_\ell)},
\end{equation}
where $\sigma_{b,y,m}=(-1)^{nm+by+\binom b2}$ is independent of $I$.

At a root $u$ of $P_z$,
\begin{equation}\label{eq:derivativeidentity0}
    uP_z'(u)=2yu^b-nz,
    \qquad n=b+y.
\end{equation}
Expanding the squared Vandermonde as
\begin{equation}\label{eq:vandsquare}
    \Delta(u_1,\ldots,u_b)^2
    =\sum_{\mathbf w}\lambda_{\mathbf w}
       u_1^{w_1}\cdots u_b^{w_b},
    \qquad
    \sum_{\ell=1}^{b}w_\ell=b(b-1),
\end{equation}
we obtain the exact decomposition
\begin{equation}\label{eq:stringdecomp}
 T_I(z)=\sigma_{b,y,m}z^m
 \sum_{\mathbf w}\lambda_{\mathbf w}
 \prod_{\ell=1}^{b}
 \frac{u_\ell^{\,1-m+w_\ell}}
      {2yu_\ell^b-nz}.
\end{equation}
This is the string decomposition used in the rest of the paper.  The weights are simply the exponents of a monomial in $\Delta(J)^2$, and \eqref{eq:vandsquare} proves the normalization
\begin{equation}\label{eq:weightsum}
    \sum_{\ell=1}^{b}w_\ell=b(b-1).
\end{equation}
The shift $1-m+w_\ell$ is likewise forced by \eqref{eq:weylcomplement}--\eqref{eq:derivativeidentity0}; it is not an additional rule.

We emphasize that \eqref{eq:stringdecomp} is an identity, not an approximation.  Summing it over the Weyl index sets of a fixed block, and then over blocks by \eqref{eq:blocksum}, reproduces $D_m$ exactly; nothing has been discarded at this stage.  The residue decomposition below therefore refines an exact identity, which is what makes \eqref{eq:combinedformula} an identity rather than merely a sufficient recipe.

After the roots in $J$ are separated according to \eqref{eq:jgfromc}, we write a weight pattern as
\begin{equation}\label{eq:string}
    r_{I_1}^{w_1}\cdots r_{I_j}^{w_j}
    R_{J_1}^{W_1}\cdots R_{J_g}^{W_g},
\end{equation}
where the labels are distinct within each root type and
\begin{equation}\label{eq:weightsum2}
    \sum_{i=1}^{j}w_i+\sum_{q=1}^{g}W_q=b(b-1).
\end{equation}
The notation records the Vandermonde weights; the actual root powers entering \eqref{eq:stringdecomp} are
\begin{equation}\label{eq:kweights}
    k_i=1-m+w_i,
    \qquad
    K_q=1-m+W_q.
\end{equation}

There is also a useful differential interpretation.  For a root branch $u(z)$ and $k\neq0$,
\begin{equation}\label{eq:derivativeidentity}
    \frac{u(z)^k}{2yu(z)^b-nz}
    =-\frac1k\frac{d}{dz}u(z)^k,
\end{equation}
which follows from $P_z'(u)u'(z)+1=0$ together with \eqref{eq:derivativeidentity0}.  The case $k=0$ is the logarithmic derivative
\begin{equation}\label{eq:logderivative}
    \frac{1}{2yu(z)^b-nz}
    =\frac{1}{u(z)P_z'(u(z))}
    =-\frac{d}{dz}\log u(z).
\end{equation}
For a large root this is again an ordinary power series.  For a small root, $\log r_I(z)=\tfrac1b\log z+O(z^{y/b})$, so \eqref{eq:logderivative} contributes a term $-1/(bz)$, that is, a factor of $z$-exponent $-1$; this is precisely the value returned by \eqref{eq:smallexponentpiece} below at $v_i=0$ and $k_i=0$, and it is the only way a negative factor exponent can arise.  In every case the denominator in \eqref{eq:stringdecomp} introduces no new root-label dependence: it differentiates the same Lagrange series and preserves its cyclic character.  This observation closes the bookkeeping between the Weyl denominator and the root-power expansions below.

\subsection{Lagrange phases and residue indices}\label{subsec:lagrangephases}

Let
\[
    \zeta_b=e^{2\pi i/b},
    \qquad
    \zeta_y=e^{2\pi i/y}.
\]
Label the small roots so that
\[
    r_I(z)\sim 2^{-1/b}\zeta_b^I z^{1/b},
    \qquad 0\le I<b,
\]
and the large roots so that
\[
    R_J(0)=2^{1/y}\zeta_y^{-J},
    \qquad 0\le J<y.
\]
Lagrange--B\"urmann inversion gives, for every fixed integer power $k$, expansions of the form
\begin{align}
 r_I(z)^k
 &=\sum_{q\ge0}A_{q,k}\,
   \zeta_b^{I(yq+k)}z^{(yq+k)/b},
 \label{eq:smallphase}\\
 R_J(z)^k
 &=B_{0,k}\zeta_y^{-Jk}
   +\sum_{q\ge1}B_{q,k}\,
   \zeta_y^{J(bq-k)}z^q,
 \label{eq:largephase}
\end{align}
where $A_{q,k}$ and $B_{q,k}$ are independent of the labels $I,J$.  Differentiating these series according to \eqref{eq:derivativeidentity} lowers every $z$-exponent by one and leaves the displayed characters unchanged.  It also annihilates the constant term $B_{0,k}\zeta_y^{-Jk}$ of \eqref{eq:largephase}; this is why only the indices $q\ge1$ contribute to a large-root factor, and why the index $V$ below is defined from $q-1$.

For a small-root factor, write $q=bp+v$ with $0\le v<b$.  For a positive-order large-root factor, write $q-1=yp+V$ with $0\le V<y$.  In both cases $p\ge0$ is the quotient and the residue $v$ or $V$ is what the character sees.  Using \eqref{eq:kweights}, the label characters become
\begin{align}
 \alpha_i&\equiv yv_i+w_i+1-m\pmod b,
 \label{eq:alpha}\\
 \beta_q&\equiv bV_q-W_q+b-1+m\pmod y.
 \label{eq:beta}
\end{align}
The associated $z$-exponents, before the common increments by multiples of $y$ coming from $p$, are
\begin{align}
 &\frac{yv_i+w_i+1-m-b}{b}
 &&\text{for a small-root factor},\label{eq:smallexponentpiece}\\
 &V_q
 &&\text{for a large-root factor}.\label{eq:largeexponentpiece}
\end{align}
Together with the global factor $z^m$ in \eqref{eq:stringdecomp}, these expressions lead to the candidate exponent defined below.

The derivation also verifies the label-independence needed for the orbit argument.  For a fixed ordered weight pattern and fixed residue data, the constants $\lambda_{\mathbf w}$, the differentiated Lagrange coefficients, and all multiplicities depend on the positions and weights but not on the numerical root labels.  Simultaneous translation of all small labels modulo $b$, or all large labels modulo $y$, merely permutes the admissible label tuples.  Hence the coefficient attached to an ordered tuple is translation invariant.

\subsection{A global residue projection}

The same root-of-unity symmetry used in \cref{lem:blockregularity} applies to numerator determinants as well.

\begin{proposition}[Shifted support]\label{prop:shiftedsupport}
Let $F(z)\in\Q[[z]]$ be a complete root expression that is homogeneous of total degree $H$ in the roots of $P_z$.  Then
\begin{equation}\label{eq:shiftedclass}
    F(z)\in z^\rho\Q[[z^y]],
    \qquad
    \rho\equiv b^{-1}H\pmod y,
    \quad 0\leq\rho<y.
\end{equation}
\end{proposition}

\begin{proof}
The scaling identity gives
\[
    F(\omega^b z)=\omega^H F(z)
\]
for every $\omega^y=1$.  Comparing coefficients yields $bk\equiv H\pmod y$ for every nonzero coefficient $[z^k]F$.  Here coprimality is used to invert $b$ modulo $y$.
\end{proof}

Equivalently, if $\eta$ is a primitive $y$th root of unity, the projection
\begin{equation}\label{eq:projection}
    \Pi_\rho F(z)
    =\frac1y\sum_{q=0}^{y-1}\eta^{-\rho q}F(\eta^qz)
\end{equation}
retains exactly the coefficients with exponent congruent to $\rho$ modulo $y$.  The denominator has $H=y(m-b)$ and hence $\rho=0$.  Numerators may occupy a shifted residue class.

\subsection{Complete character orbits}\label{subsec:orbits}

For a fixed weight pattern and residue tuple, the dependence on the root labels is the character
\begin{equation}\label{eq:character}
 \chi_{\mathbf v,\mathbf V}(\mathbf I,\mathbf J)
 =\zeta_b^{\sum_{i=1}^{j}I_i\alpha_i}
  \zeta_y^{\sum_{q=1}^{g}J_q\beta_q},
\end{equation}
with $\alpha_i,\beta_q$ as in \eqref{eq:alpha}--\eqref{eq:beta}.  By the derivation in \cref{subsec:strings,subsec:lagrangephases}, the remaining coefficients are independent of the numerical labels and are invariant under simultaneous translation.

\begin{lemma}[Translation character test]\label{lem:translation}
Let $N\geq1$, let $\zeta_N$ be a primitive $N$th root of unity, and let $c(i_1,\ldots,i_q)$ be a coefficient on ordered distinct $q$-tuples in $\Z/N\Z$ satisfying
\[
    c(i_1+1,\ldots,i_q+1)=c(i_1,\ldots,i_q).
\]
Then
\begin{equation}\label{eq:theta}
 \Theta_N(k_1,\ldots,k_q)
 =\sum_{\substack{i_1,\ldots,i_q\in\Z/N\Z\\
                   i_1,\ldots,i_q\text{ distinct}}}
 c(i_1,\ldots,i_q)
 \zeta_N^{k_1i_1+\cdots+k_qi_q}
\end{equation}
vanishes unless
\begin{equation}\label{eq:charsumzero}
    k_1+\cdots+k_q\equiv0\pmod N.
\end{equation}
\end{lemma}

\begin{proof}
Translation of every label by one preserves both the summation set and, by the preceding derivation, the coefficient $c$, while multiplying the sum by $\zeta_N^{k_1+\cdots+k_q}$.  If this multiplier is not $1$, the sum is zero.
\end{proof}

Applying the lemma independently to the small and large labels gives the modulo restrictions.

\begin{corollary}[Small- and large-label congruences]\label{cor:labelcongruences}
A complete residue orbit associated with \eqref{eq:string} can be nonzero only if
\begin{align}
 \sum_{i=1}^{j}(yv_i+w_i)+j(1-b-m)
 &\equiv0\pmod b,
 \label{eq:smallcongruence}\\
 b\sum_{q=1}^{g}V_q-\sum_{q=1}^{g}W_q
   +g(b-1+m)
 &\equiv0\pmod y.
 \label{eq:largecongruence}
\end{align}
\end{corollary}

\begin{proof}
The first relation is $\sum\alpha_i\equiv0\pmod b$; subtracting $jb$ gives the displayed form.  The second is $\sum\beta_q\equiv0\pmod y$.
\end{proof}

Define the \emph{residue offset}
\begin{equation}\label{eq:Qdef}
 \mathcal Q=
 \sum_{q=1}^{g}V_q
 +\frac1b\left[
   \sum_{i=1}^{j}(yv_i+w_i)+j(1-b-m)
 \right].
\end{equation}
When \eqref{eq:smallcongruence} holds, $\mathcal Q$ is an integer.  Using \eqref{eq:weightsum2} and $j+g=b$ gives the exact identity
\begin{equation}\label{eq:exactQidentity}
 b(\mathcal Q+m)
 =y\sum_{i=1}^{j}v_i
  +b\sum_{q=1}^{g}V_q
  -\sum_{q=1}^{g}W_q
  +g(b-1+m).
\end{equation}
Reducing \eqref{eq:exactQidentity} modulo $y$ and using the invertibility of $b$ modulo $y$ shows that \eqref{eq:largecongruence} is equivalent to
\begin{equation}\label{eq:Qcongruence}
    \mathcal Q+m\equiv0\pmod y.
\end{equation}
Thus the two congruences are precisely the integrality and character-orthogonality conditions for a residue subblock.

\begin{remark}[Why termwise rationality is not the criterion]\label{rem:rationality}
An individual phase in \eqref{eq:character} may be irrational while its complete Galois orbit has a nonzero rational sum.  For example,
\[
    \cos\frac{2\pi}{5}+\cos\frac{4\pi}{5}=-\frac12.
\]
Hence one must sum complete label orbits.  Retaining only individually rational phases can discard a legitimate rational contribution.
\end{remark}

\subsection{Candidate exponents and further cancellation}\label{subsec:candidate}

For a fixed block, weight pattern, and residue tuple satisfying the small-root integrality condition, \eqref{eq:smallexponentpiece}--\eqref{eq:largeexponentpiece} show that the corresponding complete residue contribution is supported in
\begin{equation}\label{eq:candidatesupport}
    z^{E_0}\Q[[z^y]],
    \qquad
    E_0:=m+\mathcal Q.
\end{equation}
The number $E_0$ is a \emph{candidate initial exponent}, or equivalently a lower bound on the true initial exponent after the complete orbit sum.  It need not be the actual valuation.  Write
\begin{equation}\label{eq:residueseries}
    \mathcal B_{m,c,\lambda,a}(z)
    =z^{E_0}\sum_{\ell\ge0}C_{m,c,\lambda,a}(\ell)z^{y\ell}.
\end{equation}
Then
\begin{equation}\label{eq:truevaluation}
    \ord_z\mathcal B_{m,c,\lambda,a}
    =E_0+y\min\{\ell:C_{m,c,\lambda,a}(\ell)\neq0\},
\end{equation}
whenever the subblock is nonzero.  The two congruences are necessary but not sufficient for $C(0)\neq0$: additional cancellation among admissible label permutations or weight terms may raise the valuation.  This is why the examples below, where $E_0$ happens to be attained, should not be read as a general equality statement.

If either character congruence fails, \cref{lem:translation} makes every coefficient in \eqref{eq:residueseries} zero, because increasing the Lagrange indices by full periods does not change the residue character.  If both congruences hold, the complete orbit coefficients $C(\ell)$ are evaluated only for those $\ell$ that can still lie in the finite polynomial window.

\section{The combined pruning theorem}

\begin{theorem}[Block-and-residue pruning]\label{thm:combined}
Let $D_m$ be decomposed into complete root blocks and then into residue subblocks as in \cref{sec:residues}.  For a subblock indexed by $(c,\lambda,a)$, let $E_0=m+\mathcal Q$ be its candidate exponent.  To recover $D_m$ through its degree bound $N_m$, the following steps are sufficient:
\begin{enumerate}
\item discard every root block with $L_{m,c}>N_m$;
\item within the surviving blocks, discard every residue tuple that fails either \eqref{eq:smallcongruence} or \eqref{eq:largecongruence};
\item discard every remaining residue subblock with $E_0>N_m$;
\item for each retained residue subblock, compute the complete orbit coefficients
\[
    C_{m,c,\lambda,a}(\ell),
    \qquad
    0\leq\ell\leq
    \left\lfloor\frac{N_m-E_0}{y}\right\rfloor,
\]
and discard the subblock if all of these vanish.
\end{enumerate}
Consequently,
\begin{equation}\label{eq:combinedformula}
 D_m(z)
 =[z^{\leq N_m}]
 \sum_{\substack{(c,\lambda,a)\ \mathrm{survive}\ (1)-(3)}}
 z^{E_0}
 \sum_{0\leq\ell\leq(N_m-E_0)/y}
 C_{m,c,\lambda,a}(\ell)z^{y\ell}.
\end{equation}
In particular, the first possible exponent of a surviving residue subblock is $E_0$, but its true valuation may be larger if $C(0)=0$.
\end{theorem}

\begin{proof}
By \eqref{eq:stringdecomp} and \eqref{eq:blocksum}, the sum of all residue subblocks equals $D_m$; the steps below only remove subblocks that contribute nothing in the range $0\leq k\leq N_m$.  Step (1) is \cref{cor:blockexclusion}.  Step (2) follows from \cref{lem:translation,cor:labelcongruences}; a failed character condition annihilates the complete residue orbit coefficient-by-coefficient.  Step (3) uses the lower support bound \eqref{eq:candidatesupport}.  Once $E_0\leq N_m$, the support progression in \eqref{eq:residueseries} shows that only the finitely many indices in step (4) can enter the polynomial window.  Formula \eqref{eq:combinedformula} is the resulting truncation.
\end{proof}

For a numerator of root degree $H$, the same statement holds after replacing the denominator residue class $0\pmod y$ by the shifted class $\rho\equiv b^{-1}H\pmod y$ from \cref{prop:shiftedsupport}.

The practical order of calculation is therefore
\begin{equation}\label{eq:algorithm}
\begin{gathered}
\text{whole-block valuation}\\
\downarrow\\
\text{character congruences on residue tuples}\\
\downarrow\\
\text{candidate-exponent cutoff}\\
\downarrow\\
\text{complete orbit coefficients only inside the finite window}.
\end{gathered}
\end{equation}

\section{An overlap example: \texorpdfstring{$(+3,-2)$}{(+3,-2)} with \texorpdfstring{$m=12$}{m=12}}\label{sec:example12}

Take
\[
    y=3,
    \qquad b=2,
    \qquad m=12.
\]
Then
\[
    \mathcal L=\{0,1\},\qquad
    \mathcal T_{12}=\{2,3,\ldots,11\},\qquad
    \mathcal U_{12}=\{12,13,14\}.
\]
For the hitting-time interpretation we start at $s=m-1=11$ and take boundary payoff $1$ at the inner left state $1$ and $0$ at every other absorbing state, as illustrated in \cref{fig:model}.  The normalized kernel is
\[
    P_z(X)=X^5-2X^2+z.
\]
The Weyl expansion contains
\[
    \binom52=10
\]
terms, divided into blocks of sizes
\[
    1,\quad6,\quad3.
\]
The polynomial window is
\begin{equation}\label{eq:N12}
    N_{12}=3\left\lfloor\frac{12-2}{5}\right\rfloor=6.
\end{equation}
The valuation bounds are
\begin{equation}\label{eq:L12}
    L_{12,0}=0,
    \qquad
    L_{12,1}=3\left\lceil\frac{11}{6}\right\rceil=6,
    \qquad
    L_{12,2}=3\left\lceil\frac{12}{3}\right\rceil=12.
\end{equation}
Hence the three-term all-small block is omitted, while the mixed block can meet the denominator only at the top degree $z^6$.

The global projection permits only the powers $1,z^3,z^6$ in the polynomial window.  Expanding the single all-large block only at those degrees gives
\begin{equation}\label{eq:base12}
    B_{12,0}(z)
    =1024-320z^3+3z^6+O(z^9).
\end{equation}

For the mixed block, the only possible contribution in the window is its leading residue orbit.  Here $c=1$, so by \eqref{eq:jgfromc} the string has $j=1$ small and $g=1$ large root, and the relevant weight pattern is
\[
    w=W=1,
    \qquad
    v=V=0.
\]
The character exponents are
\[
    \alpha=3v+w+1-m=-10\equiv0\pmod2,
\]
\[
    \beta=2V-W+1+m=12\equiv0\pmod3.
\]
Moreover,
\begin{equation}\label{eq:Q12}
    \mathcal Q=V+\frac12\bigl(3v+w+1-2-m\bigr)=-6.
\end{equation}
Thus the candidate exponent is $E_0=m+\mathcal Q=6$.  In this example the leading complete orbit coefficient is nonzero, so the candidate exponent is attained.  Summing the complete orbit over the two small and three large labels gives
\begin{equation}\label{eq:mixed12}
    B_{12,1}(z)=-2z^6+O(z^9).
\end{equation}
The two remaining weight patterns, $(w,W)=(2,0)$ and $(0,2)$, force $v=1$ and give candidate exponents $E_0=9$, so no other mixed residue subblock can enter the range $0\leq k\leq6$.  Therefore
\begin{equation}\label{eq:D12}
 \begin{aligned}
    D_{12}(z)
    &=[z^{\leq6}]
      \bigl(B_{12,0}(z)+B_{12,1}(z)\bigr)\\
    &=1024-320z^3+z^6.
 \end{aligned}
\end{equation}

Returning to the ordinary time variable by $z^3=t^5$, \cref{prop:transferbridge} gives
\begin{equation}\label{eq:transfer12}
    \det(I-tQ_{12})
    =1-\frac{5}{16}t^5+\frac{1}{1024}t^{10}.
\end{equation}
Solving the boundary-value problem \eqref{eq:resolvent} for the start $s=11$ and target state $1$ gives the completely self-contained hitting-time generating function
\begin{equation}\label{eq:F12physical}
    F_{11}(t)
    =\frac{t^5(t^5+32)}{t^{10}-320t^5+1024}.
\end{equation}
Its leading term is $t^5/32$, corresponding to the unique shortest path
\[
    11\to9\to7\to5\to3\to1,
\]
five consecutive $-2$ steps, each of probability $1/2$.

Thus the block bound explains why exactly one mixed orbit is needed: the base block supplies the first three admissible normalized powers, the mixed block corrects the coefficient at the unique overlap degree, and the terminal block is entirely irrelevant.

\section{A no-overlap example: \texorpdfstring{$(+5,-2)$}{(+5,-2)} with \texorpdfstring{$m=22$}{m=22}}\label{sec:example22}

Now take
\[
    y=5,
    \qquad b=2,
    \qquad m=22.
\]
Here
\[
    \mathcal L=\{0,1\},\qquad
    \mathcal T_{22}=\{2,3,\ldots,21\},\qquad
    \mathcal U_{22}=\{22,23,24,25,26\}.
\]
Again take $s=m-1=21$ and boundary payoff $1$ at state $1$ and $0$ at the other absorbing states.  The normalized kernel is
\[
    P_z(X)=X^7-2X^2+z.
\]
The Weyl expansion contains
\[
    \binom72=21
\]
terms, with block sizes
\[
    1,\quad10,\quad10.
\]
The denominator degree window is
\begin{equation}\label{eq:N22}
    N_{22}=5\left\lfloor\frac{22-2}{7}\right\rfloor=10.
\end{equation}
The two nonbase blocks begin no earlier than
\begin{equation}\label{eq:L22}
    L_{22,1}
    =5\left\lceil\frac{21}{10}\right\rceil=15,
    \qquad
    L_{22,2}
    =5\left\lceil\frac{22}{5}\right\rceil=25.
\end{equation}
Both are therefore excluded before any residue subblocks are enumerated.  Twenty of the twenty-one Weyl terms are omitted, and
\begin{equation}\label{eq:onlybase22}
    D_{22}(z)=[z^{\leq10}]B_{22,0}(z).
\end{equation}

The modulo-five projection says that only the powers $1,z^5,z^{10}$ can occur.  Expanding the single all-large term at those three powers gives
\begin{equation}\label{eq:D22}
    D_{22}(z)
    =1048576-294912z^5+5952z^{10}
    =64\bigl(16384-4608z^5+93z^{10}\bigr).
\end{equation}
No coefficients of $z^{15},z^{20}$, or higher need be computed, even though such terms occur in the separate algebraic tails and cancel in the full determinant.

Returning to the ordinary time variable by $z^5=t^7$, \cref{prop:transferbridge} gives
\begin{equation}\label{eq:transfer22}
    \det(I-tQ_{22})
    =1-\frac{9}{32}t^7+\frac{93}{16384}t^{14}.
\end{equation}
For the explicitly specified start $s=21$ and target state $1$, \eqref{eq:resolvent} gives
\begin{equation}\label{eq:F22physical}
    F_{21}(t)
    =\frac{t^{10}(5t^7+64)}
    {4(93t^{14}-4608t^7+16384)}.
\end{equation}
The leading term is $t^{10}/1024$, corresponding to ten consecutive $-2$ steps from $21$ to $1$.  Thus the probabilistic meaning of the example no longer depends on any external choice of starting position or boundary convention.

The normalized denominator occupies the residue class $0\pmod5$.  More generally, a numerator determinant may occupy a shifted residue class as described in \cref{prop:shiftedsupport}.

\section{Independent polynomial checks}

The examples can be verified without expanding any algebraic root.  Let $h_k$ be the complete homogeneous symmetric functions in the $n=b+y$ roots.  From \eqref{eq:especial},
\begin{equation}\label{eq:Hseries}
    \sum_{k\geq0}h_k u^k
    =\frac{1}{1-2u^y+zu^{b+y}},
\end{equation}
so
\begin{equation}\label{eq:hrecurrence}
    h_k=2h_{k-y}-zh_{k-b-y},
    \qquad
    h_0=1,
    \quad h_k=0\ (k<0).
\end{equation}
The Jacobi--Trudi formula gives
\begin{equation}\label{eq:JT}
    D_m(z)=
    \det\bigl(h_{m-b-i+j}\bigr)_{1\leq i,j\leq y}.
\end{equation}
For $(y,b,m)=(3,2,12)$, this yields
\[
    D_{12}(z)=1024-320z^3+z^6.
\]
For $(y,b,m)=(5,2,22)$, it yields
\[
    D_{22}(z)=1048576-294912z^5+5952z^{10}.
\]
These computations verify the final polynomials independently of the root-block construction.  The pruning method explains why only the displayed block and residue data are needed to construct them.

\section{Conclusion}

The known Schur formulas give closed rational expressions for walks between two barriers.  The contribution here is a finite construction principle for their kernel-root expansions.  The root terms are first assembled into complete small-root/large-root blocks.  A valuation comparison removes every block that begins above the degree of the desired polynomial.  Surviving blocks are then divided into residue subblocks, and complete cyclic label orbits remove all nontrivial characters.  Finally, each surviving series is truncated at the known polynomial degree.

The two examples illustrate both outcomes.  For $(+3,-2)$ with $m=12$, the mixed block reaches exactly one admissible degree, and its complete six-term orbit changes the top coefficient.  For $(+5,-2)$ with $m=22$, both nonbase blocks lie outside the polynomial window, so twenty of the twenty-one Weyl terms and all of their residue expansions are unnecessary.

This viewpoint also clarifies the role of the modulo relations in the string construction.  They are not tests of whether an individual trigonometric phase is rational.  They are character-orthogonality conditions for complete root-label orbits.  In combination with block valuation, they provide a rigorous and computationally economical pruning algorithm for two-barrier generating functions.

A natural next problem is to sharpen the valuation bound by determining when the first character-admissible coefficient of a complete block is nonzero.  Such a result would classify exact overlap degrees for broader families of step sizes while retaining the same block-and-residue framework.

\end{document}